\documentclass[review,hidelinks,onefignum,onetabnum]{siamart220329}

\usepackage{color,soul}
\usepackage{amssymb,amsmath} 
\usepackage{mathtools}
\usepackage{mathrsfs}
\usepackage{caption}
\usepackage{subcaption}

\usepackage{upgreek}
\usepackage{color,soul}
\definecolor{LightGray}{gray}{0.8}
\definecolor{Orange}{rgb}{1.0, 0.31, 0.0}
\definecolor{Green}{rgb}{0.3, 1.0, 0.3}
\definecolor{Blue}{rgb}{0.75,0.75,1}

\newcommand{\rcyl}{\rho_{\rm cyl}}
\newcommand{\raux}{\rho_{\rm aux}}
\newcommand{\rcri}{\rho_{\rm cri}}
\newcommand{\rfil}{\rho_{\rm fil}}
\newcommand{\dref}{d_{\rm ref}}

\newcommand{\xifil}{\xi_{\rm fil}}

\usepackage{lipsum}
\usepackage{amsfonts}
\usepackage{graphicx}
\usepackage{epstopdf}
\usepackage{algorithmic}
\ifpdf
  \DeclareGraphicsExtensions{.eps,.pdf,.png,.jpg}
\else
  \DeclareGraphicsExtensions{.eps}
\fi

\newsiamremark{remark}{Remark}
\newsiamremark{hypothesis}{Hypothesis}
\crefname{hypothesis}{Hypothesis}{Hypotheses}
\newsiamthm{claim}{Claim}

\headers{Convergence, oscillations in Laplace-Neumann problems}{G. D. Kolezas, G. Fikioris, and J. A. Roumeliotis}

\title{Convergence, divergence, and inherent oscillations in MFS solutions of two-dimensional Laplace-Neumann problems}

\author{Georgios D. Kolezas \thanks{School of Electrical and Computer Engineering, National Technical University of Athens (\email{geokolezas@central.ntua.gr}, \email{gfiki@ece.ntua.gr}, \email{jroum@mail.ntua.gr}).}
\and George Fikioris \footnotemark[1]
\and John A. Roumeliotis \footnotemark[1]}

\usepackage{amsopn}

\ifpdf
\hypersetup{
  pdftitle={An Example Article},
  pdfauthor={D. Doe, P. T. Frank, and J. E. Smith}
}
\fi

\begin{document}

\maketitle

\begin{abstract}
The method of fundamental solutions (MFS), also known as the method of auxiliary sources (MAS), is a well-known computational method for the solution of boundary-value problems. 
The final solution (``MAS solution'') is obtained once we have found the amplitudes of $N$ auxiliary ``MAS sources.'' Past studies have demonstrated that it is possible for the MAS solution to converge to the true solution even when the $N$ auxiliary sources diverge and oscillate. The present paper extends the past studies by demonstrating this possibility within the context of Laplace's equation with Neumann boundary conditions. One can thus obtain the correct solution from sources that, when $N$ is large, must be considered unphysical.  We carefully explain the underlying reasons for the unphysical results, distinguish from other difficulties that might concurrently arise, and point to significant differences with time-dependent problems that were studied in the past.
\end{abstract}



\section{Introduction}\label{intro}
The method of fundamental solutions (MFS) is a popular numerical technique for the solution of certain types of boundary-value problems \cite{bog_85,fai_kar_98,fai_kar_mar_03,che_hon_20}. The key concept of the method is the representation of the solution as a superposition of fundamental solutions of the governing differential equation whose singularities lie outside the problem's domain. In the context of electromagnetics, the MFS is also known as the method of auxiliary sources (MAS) and has been extensively used for the modeling of two-dimensional (2D) scattering problems \cite{kak_ana_02}. 
In its simplest form, the 2D scatterer has a perfect electric or magnetic conductivity and is illuminated externally by a point source (exterior problem), or is a cavity surrounded by a perfectly conducting electric or magnetic space, and illuminated by a source inside the cavity (interior problem). The central idea of MAS is to represent the scattered electric/magnetic field or potential by $N$ discrete sources (``MAS sources'') on a fictitious closed surface located in the complementary region. The amplitudes of the MAS sources are determined by satisfying the boundary condition on $N$ boundary points; this amounts to solving an $N\times N$ linear system. Once the system is solved, the field or potential generated by the $N$ MAS sources within the region of interest is easily determined. As with any numerical method, one expects this field or potential to converge to the true scattered field or potential as $N\to\infty$.
    
It has been known for some years \cite{fik_06, fik_psar_07, fik_tsi_15, tsi_zou_fik_lev_18,fikioris2018bookchapter}\footnote{See also \cite{mal_yer_77,bob_tom_95,bar_bet_08}, which are pioneering or independent works by other research groups. These works  are discussed in the Introduction of \cite{fikioris2018bookchapter} and in \cite{sko_13}. Ref. \cite{fikioris2018bookchapter} also discusses an interesting similarity to antenna superdirectivity, which is further elaborated in \cite{and_fik_12} and \cite{fik_pap_mav_13}} that the convergence issue is far from simple, because the MAS sources (which are intermediate results) might diverge and oscillate, even when the final result (field or potential obtained from the MAS sources) converges. Nonetheless, the phenomenon of MAS-source oscillations often goes unnoticed; this is true even for papers aiming to discuss the convergence of MFS/MAS, an example being the recent review article~\cite{che_hon_20}.

Oscillations also go unnoticed in the recent work \cite{li_lee_hua_chi_17}, which considers the stability and errors of MFS within the context of the Laplace equation with Neumann boundary conditions. Ref. \cite{li_lee_hua_chi_17} focuses mainly on the exponential growth of the condition number and the improvement of stability.  According to \cite{li_lee_hua_chi_17}, the analysis therein ``is intriguing due to its distinct features.'' 

The main purpose of the present paper  is to extend the considerations of \cite{fik_06, fik_psar_07,fik_tsi_15, tsi_zou_fik_lev_18,fikioris2018bookchapter} to the case of  Laplace's equation with Neumann boundary conditions. Unlike \cite{li_lee_hua_chi_17}, we do not focus on the condition number. We rather examine: (i) the divergence of intermediate results; (ii) the
concurrent convergence of the final solution (but, for the first time, there arises a case (\cref{internal_problem}) where the potential/final solution in a certain sense diverges); and (iii) the oscillations, which are the manifestation of the MAS-source divergence. 
We emphasize that the oscillations are unrelated to certain other errors that arise---including errors due to roundoff and matrix ill-conditioning\footnote{Let us add that, in time-harmonic cases, the oscillations are also unrelated  to the much-discussed phenomenon of internal resonances \cite{fik_06, fik_psar_07,fik_tsi_15, tsi_zou_fik_lev_18,fikioris2018bookchapter}.}---and compare to corresponding theoretical results in the important early paper of Kitagawa \cite{kit_91}. 

To the best of our knowledge, this is the first work discussing MAS-source oscillations within the specific context of Laplace-Neumann problems. As in our previous works \cite{fik_06, fik_psar_07,fik_tsi_15,fikioris2018bookchapter} we bring out most conclusions by means of an analytical study of a simple problem---with a circular boundary---in which the entire MAS solution can be found analytically, by both exact and asymptotic (large-$N$) formulas. The formulas herein, however, turn out to be simpler than those of our previous works; compare, for example, the present exact formula \eqref{AscN1} to the corresponding formula (37) of \cite{fik_06} (the latter formula pertains to the MAS electric field of a time-harmonic problem). This relative simplicity  allows us to bring out our findings in a more direct manner. Furthermore, we find that the present (Laplace-Neumann) case presents several features that do not arise in the time-harmonic ones; one example is the aforementioned ``divergence'' of the MAS potential.

We use the language of magnetostatics for definiteness, but our Neumann problems are generic. We specifically deal with 2D exterior (or interior) Neumann problems satisfied by the vector potential $\mathbf{A}={\mathbf{\hat z}}A_z$ within unbounded (or bounded) regions of permeability $\mu=\mu_0=4\pi\times 10^{-7}$~H/m (SI units are used throughout), with the complementary bounded (or unbounded) region having infinite permeability ($\mu=\infty$). The source is an infinitely long, temporally and spatially constant electric current $I$ (units: A). Our most detailed results (sections~\ref{section: preliminaries}--\ref{section:numerics}) concern the exterior circular problem: we consider two solutions to this problem, and compare the two in detail. In \cref{section-extensions}, we extend to other Laplace-Neumann problems, including a noncircular problem with an elliptic boundary. As one might expect the main difficulty---i.e., the oscillations---that occurs in the simple circular problem is also present in the more complicated one.  

Now let us, for definiteness, refer to exterior problems, so that the auxiliary surface is placed in the interior. Exterior potentials generated by smooth, continuous, 2D surface current densities $K_z$ (units: A/m) are analytic functions of the space variables. Therefore no such $K_z$ can generate an $A_z$ exhibiting a singularity at a point exterior to $K_z$. MAS-source oscillations occur when the aforementioned auxiliary surface lies behind a singularity of the analytic continuation of the scattered potential (the analytic continuation is performed into the scatterer's interior): When this is the case, the $N$ discrete currents---when properly normalized---cannot converge to a $K_z$. While this is easy to understand, it is less apparent that \textit{divergent} currents (even when normalized) are capable of generating a \textit{convergent} potential.  In most realistic problems, the singularities of the analytic continuation are not known beforehand. A circular boundary has the advantage of bringing out the aforesaid convergence unequivocally, because there are image singularities, and they are known explicitly. Thus much herein pertains to the question: What happens when the auxiliary surface \textit{does not} enclose all singularities?, and to the corresponding question for interior problems.

In all our Neumann problems, the so-called self-consistency condition (i.e., the condition that must be satisfied by the integral of the Neumann boundary data for the problem to have a solution see, e.g., \cite{barton}) can be viewed as a consequence of Amper\'e's circuital law, which states that the line integral of the $\mathbf{H}$-field around a closed loop equals the total current passing through the loop. For clarity, however, we verify all self-consistency conditions by independent means. 

Throughout this paper, the terms MAS and MFS are used interchangeably. 

\section{Preliminaries}
\label{section: preliminaries}
We begin by giving some interrelated series and integrals to be used throughout. First consider the convergent definite integral defined by
\begin{multline}\label{idefinition}
J(m,x,y)=\int_{-\pi}^{\pi}\frac{(1-x\cos\theta)\cos(m\theta)}{y^2+(1-y^2)\cos^2\theta-2x\cos\theta +x^2}{\rm d}\theta, \\ \quad 0<y\leq1<x,\quad m=0,1,2,\ldots.
\end{multline}
The integral equals
\begin{align}
\label{ievaluated}
J(m,x,y)=\frac{\pi\left[(1-y)^m-(1+y)^m\right]}{y\left(x+\sqrt{x^2+y^2-1}\right)^m},
\end{align}
but we were not able to find this evaluated form in the usual sources \cite{prudnikov, gradshteyn}. We thus prove \eqref{ievaluated} in \cref{appendix-integral}. Our next two integrals are corollaries involving the Poisson kernel, 
\begin{align}
\label{ievaluated-circle-m-equals-zero}
\int_{-\pi}^{\pi}\frac{1-x\cos\theta}{x^2-2x\cos\theta +1}{\rm d}\theta=\begin{cases}
0,\qquad &x>1\\
2\pi,\qquad 0\le\!\!\!\!&x<1;
\end{cases}
\end{align}
\begin{align}
\label{ievaluated-circle-m-notequal-zero}
\int_{-\pi}^{\pi}\frac{(1-x\cos\theta)\cos(m\theta)}{x^2-2x\cos\theta +1}{\rm d}\theta=\begin{cases}
-\pi x^{-m},\qquad &x>1,\quad m=1,2,\ldots\\
\pi x^m,\qquad 0\le\!\!\!\!&x<1,\quad m=1,2,\ldots.
\end{cases}
\end{align}
The top integrals (\ref{ievaluated-circle-m-equals-zero}) and (\ref{ievaluated-circle-m-notequal-zero}) are $J(m,x,1)$. The bottom ones follow immediately, because (\ref{idefinition}) gives $J(m,x,1)+J(m,1/x,1)=2\pi$ or $0$, depending on whether $m=0$ or $m\ge 1$. By setting $x=\rho_2/\rho_1$ in (\ref{ievaluated-circle-m-equals-zero}) and (\ref{ievaluated-circle-m-notequal-zero})---or by  direct summation---we can deduce the two Fourier series
\begin{align}\label{identity1}
 \frac{\rho_1-\rho_2\cos\theta}{\rho_1^2+\rho_2^2-2\rho_1\rho_2\cos\theta}=-\frac{1}{2\rho_1}\sum_{n\neq 0}\bigg(\frac{\rho_1}{\rho_2}\bigg)^{|n|}e^{in\theta}, \quad 0\le\rho_1<\rho_2;
 \end{align}
\begin{align}\label{identity2}
 \frac{\rho_1-\rho_2\cos\theta}{\rho_1^2+\rho_2^2-2\rho_1\rho_2\cos\theta}=\frac{1}{\rho_1}+\frac{1}{2\rho_1}\sum_{n\neq 0}\bigg(\frac{\rho_2}{\rho_1}\bigg)^{|n|}e^{in\theta}, \quad 0\le \rho_2<\rho_1.
 \end{align}
(Throughout this paper, $\sum_{n\neq 0}$ denotes $\sum_{n=-\infty}^{-1}+\sum_{n=1}^{\infty}$.) A third useful Fourier series is
\begin{multline}\label{identity3}
\ln{\frac{\sqrt{\rho_1^2+\rho_2^2-2\rho_1\rho_2\cos\theta}}{\rho_3}}=\ln{\frac{\rho_2}{\rho_3}}-\frac{1}{2}\sum_{n\neq 0}\frac{1}{|n|}\bigg(\frac{\rho_1}{\rho_2}\bigg)^{|n|}e^{in\theta}, \\ \quad 0\le\rho_1<\rho_2, \quad 0<\rho_3,
\end{multline}
whose $\rho_1$-derivative is tantamount to (\ref{identity1}).

We solve $N\times N$ symmetric circulant  systems by means of the Discrete Fourier transform (DFT). The DFT of  $\alpha_0,\alpha_1,\ldots,\alpha_{N-1}$ is  denoted by $\alpha^{(0)},\alpha^{(1)},\ldots,\alpha^{(N-1)}$. The transform relations are 
\begin{align}
\label{eq:DFT-relations}
    \alpha^{(m)}=\frac{1}{N}\sum_{p=0}^{N-1}\alpha_pe^{-i2\pi mp/N},\qquad
     \alpha_p=\sum_{m=0}^{N-1}\alpha^{(m)} e^{i2\pi mp/N}.
\end{align}
In (\ref{eq:DFT-relations}), the second (inverse-transform) relation follows from the first via the identity
\begin{align}\label{sumexp}
\sum_{p=0}^{N-1}e^{i2\pi mp/N}=\begin{cases}
	N,\quad \text{if}\,\,m=\text{multiple of}\,\,\,N,\\	0,\quad\text{otherwise}.
\end{cases}
\end{align}
A symmetric circulant system is one whose left-hand side is expressed as the convolution sum $\sum_{l=0}^{N-1}\alpha_{l-p}\,\beta_l$, in which $\alpha_p=\alpha_{-p}=\alpha_{p+N}$. By the definition in \eqref{eq:DFT-relations}, the DFT of this sum equals $N\alpha^{(m)}\beta^{(m)}$.

\section{Exterior circular problem: Problem formulation}
\label{Problemdesciption}
\begin{figure}
\centering
\includegraphics[scale=0.47]{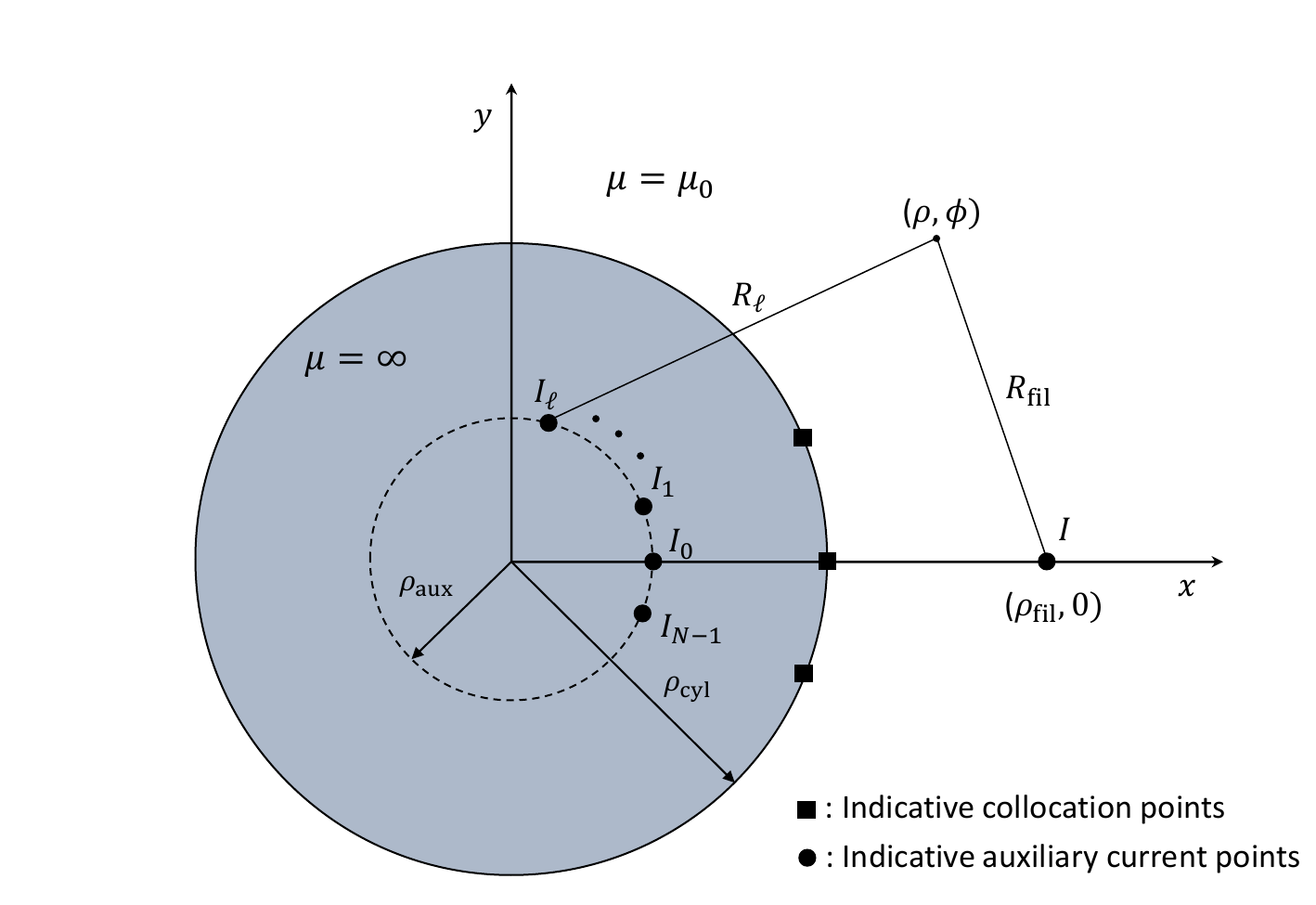}
\caption{Geometry of the exterior circular problem. Indicative positions of MAS currents and collocation points are depicted.}
\label{geometry}
\end{figure}

We consider the 2D magnetostatic problem shown in \cref{geometry}, where a circular cylinder of infinite permeability ($\mu=\infty$) is illuminated by a temporally constant line-current $I$. 
The cylinder has radius $\rcyl$, infinite extent along the $z$-direction, and axis the $z$-axis. The infinitely long current $I$ is parallel to the $z$-axis and positioned at $(\rho,\phi)=(\rho_{\rm fil},0)$, with $\rfil>\rcyl$. The medium surrounding the cylinder is free space with permeability $\mu_0$.

Because the examined geometry is invariant along the $z$-direction, the total magnetostatic vector potential $\mathbf{A}^{\rm tot}$ only has a $z$-component and depends on the polar coordinates $(\rho,\phi)$, i.e., $\mathbf{A}^{\rm tot}=\mathbf{\hat z}A_z^{\rm tot}(\rho,\phi)$. For $\rho>\rcyl$ and $(\rho,\phi)\neq(\rfil,0)$, $A_z^{\rm tot}$ satisfies Laplace's equation. Adopting terms used in electromagnetic scattering, the exterior to the cylinder $A_z^{\rm tot}(\rho,\phi)$ is the sum of the incident vector potential $A_z^{\rm inc}(\rho,\phi)$ due to the current source and the vector potential $A_z^{\rm sc}(\rho,\phi)$ scattered by the cylinder. The former is given by
\begin{align}\label{Az_inc}
A_z^{\rm inc}(\rho,\phi)=-\frac{\mu_0 I}{2\pi}\ln\frac{R_{\rm fil}}{d_{\rm ref}},
\end{align}
where 
\begin{align}\label{Rfil}
R_{\rm fil}=\sqrt{\rho^2+\rfil^2-2\rho\rfil\cos\phi} 
\end{align}
is the distance between the observation point $(\rho,\phi)$ and the current source (see \cref{geometry}), while $\dref$ is an arbitrary reference distance  (so that the ratio $R_{\rm fil}/\dref$ be dimensionless and $A_z^{\rm inc}=0$ whenever $R_{\rm fil}=\dref$).
\par
Inside the infinitely permeable cylinder the magnetic field $\mathbf{H}$ is zero and the magnetic flux density $\mathbf{B}$ is finite. This means that the tangential component of $\mathbf{H}$ is continuous across the surface of the cylinder and thus zero at $\rho=\rcyl+0$. Since $\mathbf{H}=(1/\mu_0)\nabla\times\mathbf{A}$, the tangential component of the magnetic field is $H_\phi=(-1/\mu_0)\partial A_z/\partial\rho$; therefore, the Neumann boundary condition $\partial A_z^{\rm tot}/\partial\rho=0$ must hold at $\rho=\rcyl+0$. 

Our problem can be reformulated in terms of $A_z^{\rm sc}$ as follows.  We require $A_z^{\rm sc}\rightarrow0$ as $\rho\rightarrow\infty$, consistent with the fact that no current flows inside the cylinder or on its surface (at infinity, $A_z^{\rm tot}$ behaves as does $A_z^{\rm inc}$). For $\rho>\rcyl$, $A_z^{\rm sc}$ satisfies Laplace's equation subject to the nonhomogeneous Neumann boundary condition
\begin{align}\label{bc0}
    \frac{\partial A_z^{\rm sc}}{\partial\rho}=-\frac{\partial A_z^{\rm inc}}{\partial\rho}, \,\,\, \text{at}\,\,\, \rho=\rcyl+0,
\end{align}
in which the right-hand side is known from \eqref{Az_inc} and \eqref{Rfil}. Within the loop of radius $\rho_\mathrm{cyl}$, the incident $\mathbf{H}$-field surrounds no current. Therefore the self-consistency condition for our Neumann problem (see Introduction or \cite{barton}) is
\begin{align}\label{scc_external2}
    \int_0^{2\pi}\left(-\frac{\partial A_z^{\rm inc}}{\partial\rho}\right)\,\rcyl{\rm \,d}\phi =0,\,\,\, \text{at}\,\,\, \rho=\rcyl+0.
\end{align}
We can independently verify (\ref{scc_external2}) by substituting \eqref{Az_inc} and \eqref{Rfil} into the integral of its left-hand side: The resulting integral can be evaluated using the top integral \eqref{ievaluated-circle-m-equals-zero}, and turns out to equal zero.

The boundary-value problem we formulated has a well-known solution, viz., 
\begin{align}\label{Asc-exactsolution}
A_z^{\rm sc}(\rho,\phi)= -\frac{\mu_0I}{2\pi}\ln\frac{\sqrt{\rho^2+\rcri^2-2\rho\rcri\cos\phi}}{\rho},
\end{align}
where
\begin{equation}
\label{rho-cri}
   \rcri=\frac{\rcyl^2}{\rfil}\quad \left(\rcri<\rcyl\right).
\end{equation}
\Cref{Asc-exactsolution} can be found using the method of images---there are two images, one located on the $z$-axis, and the other at $(\rho,\phi)=(\rcri,0)$---or by solving the Laplace-Neumann problem using elementary methods such as the separation of variables. 

Proceeding as if this solution were not known, the next two sections apply  two different MFS schemes to the above-described boundary-value problem. 
In both schemes, the fundamental solutions originate from sources $I_0,I_1,\ldots,I_{N-1}$ placed on an interior cylindrical surface $\rho=\rho_{\rm aux}$, so that
\begin{equation}
    \label{distance-inequality}
    \raux<\rcyl<\rfil,
\end{equation}
see \cref{geometry}. The $N$ coefficients $I_0,I_1,\ldots,I_{N-1}$ are found by satisfying the Neumann boundary condition \eqref{bc0} at $N$ collocation points on  $\rho=\rcyl$. The polar angles of both the sources and the collocation points are $\phi=2\pi p/N$, where $p=0,1,\dots,N-1$.

As we will see, both our schemes arrive at the true potential \eqref{Asc-exactsolution} as $N\to\infty$, but not always in a straightforward manner. In particular, our discussions will reveal important differences between the two schemes, as well as the crucial role of the parameter $\raux$.

\section{Exterior circular problem: MAS with bounded fundamental solutions}\label{MASkats}
The two MAS schemes differ in the choice of fundamental solutions. The first version takes 
\begin{align}\label{Asc1}
    A_{z,N}^{\rm sc}(\rho,\phi)=-\frac{\mu_0}{2\pi}\sum_{\ell=0}^{N-1}I_\ell\ln\frac{R_\ell}{\rho},
\end{align}
where 
\begin{align}\label{Rl}
R_\ell=\sqrt{\rho^2+\raux^2-2\rho\raux\cos\bigg(\phi-\frac{2\pi\ell}{N}\bigg)}  
\end{align}
is the distance between the observation point $(\rho,\phi)$ and the $\ell$-th auxiliary source $I_\ell$, as shown in \cref{geometry}. The subscript $N$ in \eqref{Asc1} stresses the fact that the MAS vector potential depends on the number of auxiliary sources. In \eqref{Asc1}, we have chosen the so-called bounded fundamental solutions $\ln (R_\ell/\rho)$, which vanish at infinity. This choice, which is based on suggestions in \cite{fik_tsi_15}, \cite{kat_89}, and \cite{li_lu_xie_12}, has the desirable feature that $A_{z,N}^{\rm sc}\to 0$ as $\rho\to\infty$, for all finite $N$. 

\subsection{Exact equations for finite $N$}
\label{set-up-katsurada}
Since
\begin{align}\label{Asc2}
    A_{z,N}^{\rm sc}(\rho,\phi)=-\frac{\mu_0}{2\pi}\sum_{\ell=0}^{N-1}I_\ell\ln\frac{R_\ell}{\dref}+\frac{\mu_0}{2\pi}\sum_{\ell=0}^{N-1}I_\ell\ln\frac{\rho}{\dref},
\end{align}
our scheme amounts to using the MAS currents $I_\ell$ with the so-called ``traditional fundamental solutions'' $\ln (R_\ell/d_\mathrm{ref})$ (more on these below), while placing an extra auxiliary current $I_N=-\sum_{\ell=0}^{N-1}I_\ell$ at the origin. $I_N$ is not an independent source but its value follows from the rest of the $I_\ell$. Note that the total sum of all the auxiliary currents (including $I_N$) is equal to zero, in compliance with the requirement that no current flows on or in the cylinder.
\par
We require $A_{z,N}^{\rm sc}$ to satisfy the boundary condition \cref{bc0} at the $N$ equispaced collocation points $(\rcyl,2\pi p/N)$, where $p=0,\dots,N-1$. Using \cref{Az_inc}, \cref{Rfil}, \cref{Rl}, \cref{Asc2}, and some algebra, we readily obtain the $N\times N$ system
\begin{equation}\label{sys1}
 \sum_{\ell=0}^{N-1}B_{\ell-p}I_\ell=D_p, \quad p=0,\dots,N-1,
\end{equation}
in which
\begin{equation} \label{Bl}
  B_{p}=\frac{\rcyl^2-\rcyl\raux\cos\frac{2\pi p}{N}}{\rcyl^2+\raux^2-2\rcyl\raux\cos\frac{2\pi p}{N}}-1,
\end{equation}
\begin{equation} \label{Dp}
    D_p=-I\frac{\rcyl^2-\rcyl\rfil\cos\frac{2\pi p}{N}}{\rcyl^2+\rfil^2-2\rcyl\rfil\cos\frac{2\pi p}{N}}.
\end{equation}
Note that $B_p$ and $D_p$ are independent of $\dref$. The system \cref{sys1}--\cref{Dp} is symmetric circulant, as defined in our Introduction, and is solved by means of the DFT in \cref{appendix-circulant}. Because the DFTs $B^{(m)}$ and $D^{(m)}$ are uncomplicated, e.g.,
\begin{align}\label{Bm_exact}
B^{(0)}=\frac{\Big(\frac{\raux}{\rcyl}\Big)^N}{1-\Big(\frac{\raux}{\rcyl}\Big)^N}, \quad B^{(m)}=\frac{1}{2}\frac{\Big(\frac{\raux}{\rcyl}\Big)^{N-m}+\Big(\frac{\raux}{\rcyl}\Big)^m}{1-\Big(\frac{\raux}{\rcyl}\Big)^N},\quad m=1,\dots,N-1, 
\end{align}
\cref{appendix-circulant} further obtains a solution that has an unusually simple form; compare, for example, to the analogous solution in \cite{fik_06} (which pertains to the electrodynamics case). The final formulas thus obtained are
\begin{align}\label{Iell}
    I_\ell=\sum_{m=0}^{N-1}I^{(m)}e^{i2\pi m\ell /N}, \quad \ell=0,1,\dots,N-1,
\end{align}
in which $I^{(m)}$, which is the DFT of the solution, is given by
\begin{align}
I^{(0)}&=\frac{I}{N}\frac{\Big(\frac{\rcyl}{\rfil}\Big)^{N}\bigg[1-\Big(\frac{\raux}{\rcyl}\Big)^N\bigg]}{\Big(\frac{\raux}{\rcyl}\Big)^{N}\bigg[1-\Big(\frac{\rcyl}{\rfil}\Big)^N\bigg]},\label{I0_exact} \\
I^{(m)}&=\frac{I}{N}\frac{\bigg[\Big(\frac{\rcyl}{\rfil}\Big)^{N-m}+\Big(\frac{\rcyl}{\rfil}\Big)^m\bigg]\bigg[1-\Big(\frac{\raux}{\rcyl}\Big)^N\bigg]}{\bigg[\Big(\frac{\raux}{\rcyl}\Big)^{N-m}+\Big(\frac{\raux}{\rcyl}\Big)^m\bigg]\bigg[1-\Big(\frac{\rcyl}{\rfil}\Big)^N\bigg]}, \quad m=1,\ldots,N-1. \label{Im_exact}
\end{align}
Since $I^{(m)}=I^{(N-m)}$ by \cref{Im_exact}, when $N=$ odd we can write \eqref{Iell} as
\begin{align}\label{Iell-cosines}
    I_\ell=I^{(0)}+2\sum_{m=1}^{(N-1)/2}I^{(m)}\cos\frac{2\pi m\ell}{N}, \quad \ell=0,1,\dots,N-1\quad (N=\mathrm{odd}),
\end{align}
in which all quantities are real. (A similar relation holds when $N=$ even.) To summarize, \cref{Iell}--\cref{Iell-cosines} give the \textit{exact} MAS currents $I_\ell$, for any finite $N$. Using these $I_\ell$ and \cref{Asc1}, we can find the associated MAS potential $A_{z,N}^{\rm sc}$. We now discuss some properties of $I_\ell$ and $A_{z,N}^{\rm sc}$ for the case where $N$ is large.

\subsection{Large-$N$ behavior of MAS currents}
\label{asymptotic}
This section assumes $N=$ odd, so that \cref{Iell-cosines} applies. As $N\to\infty$, we can use \cref{distance-inequality}, \cref{I0_exact}, and \cref{Im_exact} to obtain the asymptotic relations 
\begin{align}
I^{(0)}&\sim\frac{I}{N}\bigg(\frac{\rcri}{\raux}\bigg)^N,\label{asyI0}\\
I^{(m)}&\sim\frac{I}{N}\bigg(\frac{\rcri}{\raux}\bigg)^{m},\quad m=1,2,\ldots,\frac{N-1}{2},\label{asyIm}
\end{align}
where $\rcri$ is the critical radius already defined in \eqref{rho-cri}. We further define the important dimensionless parameter $t$ as
\begin{equation}
\label{parameter-t}
    t=\frac{\rcri}{\raux}.
\end{equation}
Substitution into \eqref{Iell-cosines} and exact evaluation of the resulting sum yields
\begin{multline}\label{identity4}
    I_\ell\sim\frac{I}{N}\left[t^N+2\sum_{m=1}^{(N-1)/2}t^m\cos\frac{2\pi m \ell}{N}\right] \\
    =\frac{I}{N}\left[t^N+2\frac{(-1)^{\ell}t^{\frac{N+1}{2}}(t-1)\cos\frac{\pi\ell}{N}-t^2+t\cos\frac{2\pi\ell}{N}}{t^2-2t\cos\frac{2\pi\ell}{N}+1}\right], \quad \text{as}\,\,\, N\rightarrow\infty.
\end{multline}
\Cref{identity4}, which is a closed-form asymptotic expression for the MAS currents, is better understood by simplifying its right-hand side. To do this we distinguish three cases as follows.
\begin{itemize}
    \item{Case 1: Physical case $\raux>\rcri$}
    
        When $t<1$ or $\rcri<\raux$, we neglect the term $t^N$, as well as the \textit{first} term in the numerator (both these terms are exponentially small). Consequently,
\begin{equation}\label{asyILodd-tsmall}
I_{\ell}\sim\frac{2I}{N}\frac{-t^2+t\cos\frac{2\pi\ell}{N}}{t^2-2t\cos\frac{2\pi\ell}{N}+1}, \quad \text{as}\,\,\, N\rightarrow\infty\quad (\raux>\rcri).
\end{equation}
As $N$ grows, the arc-length distance between the adjacent $I_{\ell}$ becomes smaller and, when normalized to $2\pi\raux/N$, the discrete  $I_{\ell}$ converge to a continuous $z$-directed surface current density $K_z(\phi)$, where $\phi=2\pi\ell/N$. This $K_z(\phi)$ is located at $\rho=\raux$. 

\item{Case 2: Unphysical/oscillating case $\raux<\rcri$}
  
In stark contrast, when $t>1$ or $\raux<\rcri$,  $t^N$ is the \textit{dominant} term, and we neglect the \textit{last two} terms in the numerator. We thus obtain
\begin{equation}\label{asyILodd}
I_{\ell}\sim\frac{I}{N}t^N+\frac{2I}{N}\frac{(-1)^{\ell}t^{\frac{N+1}{2}}(t-1)\cos\frac{\pi\ell}{N}}{t^2-2t\cos\frac{2\pi\ell}{N}+1}, \quad \text{as}\,\,\, N\rightarrow\infty\quad (\raux<\rcri).
\end{equation}
In \eqref{asyILodd}, the $(-1)^\ell$ means that the $I_\ell$ oscillate, very rapidly, about the mean value (dominant term) $\frac{I}{N}t^N$, which is exponentially large in $N$. The $t^\frac{N+1}{2}$ shows that the deviations from the mean are themselves exponentially large. Such currents (however normalized) cannot converge to a continuous $K_z(\phi)$. For this reason, the case  $\raux<\rcri$ is termed ``unphysical'' or ``oscillating.''

\item{Case 3: Critical case $\raux=\rcri$}

When $t=1$ or $\raux=\rcri$, an evaluation of the right-hand side of \eqref{identity4} (which is indeterminate when $\ell=0$ and equal to zero otherwise) yields $I_0\sim I$ and, when $\ell\ne 0$, $I_\ell\to 0$. The meaning of this large-$N$ result is clear and was to be expected, because the 
exact solution
\eqref{Asc-exactsolution} coincides with a single term of the sum in \eqref{Asc1}, namely the term with $\ell=0$ \Big(this is true because $\raux=\rcri$ implies $R_0=\sqrt{\rho^2+\rcri^2-2\rho\rcri\cos\phi}$\,\Big). Consequently, the aforementioned large-$N$ result amounts to 
\begin{equation}
\label{obtain-exact-solution}
    A_{z,N}^{\rm sc}\sim A_z^{\rm sc}, \quad \text{as}\,\,\, N\rightarrow\infty, 
\end{equation} 
with only the $\ell=0$ term surviving in the limit.
\end{itemize}
\subsection{Convergence of vector potential}
\label{convMAS1}
Having shown \cref{obtain-exact-solution} for Case~3,
we discuss the large-$N$ behavior of $A_{z,N}^{\rm sc}$ for Cases 1 and 2. Substitution of \cref{Rl} into \cref{Asc1} and use of  the identity \cref{identity3} gives
\begin{multline}\label{AscN2}
A_{z,N}^{\rm sc}(\rho,\phi)=\frac{\mu_0}{4\pi}\sum_{\ell=0}^{N-1}I_\ell\Bigg[\sum_{m\neq 0}\frac{1}{|m|}\bigg(\frac{\raux}{\rho}\bigg)^{|m|}e^{im(\phi-2\pi\ell/N)}\Bigg] \\ =\frac{\mu_0}{4\pi}\sum_{m\neq 0}\frac{1}{|m|}\bigg(\frac{\raux}{\rho}\bigg)^{|m|}e^{im\phi}\sum_{\ell=0}^{N-1}I_\ell e^{-i2\pi m\ell/N}.
\end{multline}
In the last expression, the finite sum (over $\ell$) is $N$-periodic in $m$; since $I_\ell=I_{N-\ell}$, this sum is also even in $m$. By the inverse-DFT relation in \cref{eq:DFT-relations}, the sum equals  the $N$-periodic and even quantity $NI^{(m)}$. We thus obtain the exact formula
\begin{align}\label{AscN1}
A_{z,N}^{\rm sc}(\rho,\phi)=\frac{\mu_0}{4\pi}\sum_{m\neq 0}\frac{NI^{(m)}}{|m|}\bigg(\frac{\raux}{\rho}\bigg)^{|m|}e^{im\phi},
\end{align}
in which $I^{(0)}$ does not appear. Note that the aforementioned periodicity (together with the inequality $\raux<\rho$) ensure rapid convergence of the series in \cref{AscN1} for any finite $N$, even in Case~2: The product $I^{(m)}(\raux/\rho)^{|m|}$ is always bounded by a constant times the exponentially decreasing function $(\raux/\rho)^{|m|}$.

As $N\rightarrow\infty$, the asymptotic formula \cref{asyIm} (for $I^{(m)}$) holds for  
all $m>0$ and, upon replacing $m$ by $|m|$, can  be extended to all $m\ne 0$. Consequently, \cref{AscN1} yields
\begin{align}\label{Asclim1}
A_{z,N}^{\rm sc}(\rho,\phi)\sim\frac{\mu_0I}{4\pi}\sum_{m\neq 0}\frac{1}{|m|}\bigg(\frac{\rcri}{\rho}\bigg)^{|m|}e^{im\phi},\quad \mathrm{as\ } N\to\infty.
\end{align}
The series converges because $\rcri<\rcyl<\rho$, see \cref{rho-cri}, and can be evaluated (exactly) with the aid of \cref{identity3}. The result is
\begin{align}\label{Asclim3}
A_{z,N}^{\rm sc}(\rho,\phi)\sim -\frac{\mu_0I}{2\pi}\ln\frac{\sqrt{\rho^2+\rcri^2-2\rho\rcri\cos\phi}}{\rho},\quad \mathrm{as\ } N\to\infty,
\end{align}
whose right-hand side coincides with the exact solution in \cref{Asc-exactsolution}, demonstrating that \cref{obtain-exact-solution} also holds in Cases~1 and~2.
The salient point is that one obtains the true potential $A_z^{\rm sc}$ when $\raux<\rcri$. This is peculiar, because \textit{the correct solution is arrived at even in the unphysical case, in which the MAS currents diverge and oscillate rapidly.}  

Having thus clarified the particularities of the ``unphysical'' case, let us stress that some sort of anomalous behavior was to be expected, if only from the general considerations of our Introduction.  The point $(\rho,\phi)=(\rcri,0)$ is a singularity of the analytic continuation of the (exact) scattered potential to points interior to the circle, i.e., to points $(\rho,\phi)$ with $\rho<\rcyl$. The origin is another singularity. In this simple problem, the singularity positions can be deduced a priori from \eqref{Asclim3}, and are consistent with the annulus of convergence of the Laurent series in \eqref{Asclim1}. 

\subsection{Matrix condition number}
\label{section-condition-number}

We now briefly discuss the condition number of the matrix in \eqref{sys1}. For similar results, see \cite{li_lee_hua_chi_17} for the interior circular Laplace-Neumann problem;    \cite{kit_88,kit_91,smy_kar_01} for the interior circular Laplace-Dirichlet problem;  \cite{li_hua_hua_10} for a stability analysis for the circular/noncircular interior Laplace problem with a mixed type of boundary conditions; \cite{ana_lum_kak_04} for the electrodynamic problem of scattering from a PEC cylinder, and \cite{ana_kak_04} for a similar problem with a dielectric cylinder.

The DFT $B^{(m)}$ of the first row of the circulant matrix in \eqref{sys1} is given in \eqref{Bm_exact}. As $N\to\infty$, the smallest $B^{(m)}$ is $B^{(0)}\sim (\raux/\rcyl)^N$, while the largest is $B^{(1)}\sim \frac{1}{2}\raux/\rcyl$. Since circulant matrices are normal, the $B^{(m)}$ coincide with the matrix eigenvalues and singular values \cite{davis,golub}, and the ratio $B^{(1)}/B^{(0)}$ is the 2-norm condition number $\kappa$:
\begin{equation}
\label{condition-number-kappa}
    \kappa\sim \frac{1}{2}\left(\frac{\rcyl}{\raux}\right)^{N-1},\quad \mathrm{as\ }N\to\infty.
\end{equation}
In both cases, therefore (physical and unphysical), $\kappa$ is exponentially large in $N$. The base (which, as expected, is independent of $\rfil$) is seen to decrease with $\raux$. 

\subsection{Insensitivity of vector potential}
\label{section-insensitivity}

Since $\kappa$ is large, small errors (e.g., in the calculation of the boundary data/incident vector potential, or in the potentials generated by the auxiliary sources) will greatly affect the $I_\ell$. Will these resulting errors in $I_\ell$ affect the calculation of $A_{z,N}^{\rm sc}$ (which we find from \eqref{Asc1}, and is the final quantity of interest)? This is a separate issue (see, e.g., \cite{che_hon_20})  that seems to have been first addressed quantitatively by Kitagawa in the 1991 paper \cite{kit_91} (in the context of interior Laplace-Dirichlet problems, but much carries over to the present case). Since errors are not our main focus---we rather focus on phenomena that would occur even with perfect hardware and software---we limit ourselves to some qualitative remarks for the unphysical/oscillating case (Case 2, $\raux<\rcri$). We assume,  of course, that $N$ is large.

As we already saw, the main contribution to the sum in \eqref{Iell-cosines} (for $I_\ell-I^{(0)}$) comes from large values of $m$; hence the oscillations. By contrast, the main contribution to the rapidly-converging, \textit{exact} infinite series \eqref{AscN1} (for $A_{z,N}^{\rm sc}$) comes from small values of $|m|$. This allows us to view the factor $(\raux/\rho)^{|m|}$ as a low-pass filter that  blocks the large-$|m|$ values of $I^{(m)}$. 

Let us now assume that the $N\times N$ system \eqref{sys1} is solved by a method that mostly introduces errors in the large-$|m|$ values of  $I^{(m)}$. A similar hypothesis is further discussed by  Kitagawa in \cite{kit_91}, where it is phrased in terms of the matrix singular values $B^{(m)}$.  In any case, the above assumption seems to be true for the standard MATLAB solver that we used. Subject to our assumption, the aforementioned small errors we started out with (which, due to the large $\kappa$, appear as large errors in the $I_\ell$) will be filtered out, and will not greatly affect the final result $A_{z,N}^{\rm sc}$. We call this phenomenon \textit{insensitivity}. In fact, it is easily understood from \eqref{AscN1} that more insensitivity occurs for larger $\rho$, in which case the cutoff of the low-pass filter becomes smaller. 

Let us emphasize that insensitivity pertains to errors that arise from, e.g., roundoff, or transcendental-function computations, and that appear magnified in the $I_\ell$ because of the large $\kappa$. The unphysical oscillations of $I_\ell$ (Section~\ref{asymptotic}) are definitely not errors of this sort. In a sense, the two effects compete: Since large intermediate results should in
general be avoided \cite{higham}, the existence of oscillations is a hindrance to numerically obtaining the correct potential; by contrast, insensitivity helps.

\section{Exterior circular problem: Traditional fundamental solutions}\label{MAStrad}

We now employ MAS with the ``traditional'' or ``logarithmic'' fundamental solutions \cite{fik_tsi_15,kat_89,li_lu_xie_12}. This means approximating the scattered vector potential by
\begin{align}\label{Asc3}
    A'^{\rm sc}_{z,N}(\rho,\phi)=-\frac{\mu_0}{2\pi}\sum_{\ell=0}^{N-1}I'_\ell\ln\frac{R_\ell}{\dref},
\end{align}
with the $R_\ell$ defined in \eqref{Rl}. Once again, we require the $I'_\ell$ to satisfy the boundary condition \eqref{bc0} at $(\rcyl,2\pi p/N)$. It is apparent from \eqref{Asc2} that $I'_\ell\ne I_\ell$. Following the steps of Section \ref{set-up-katsurada}, we obtain a new linear system of equations that reads
\begin{equation}\label{sys2}
 \sum_{\ell=0}^{N-1}B'_{\ell-p}I'_\ell=D_p, \quad p=0,\dots,N-1,
\end{equation}
with
\begin{equation} \label{}
  B'_p=\frac{\rcyl^2-\rcyl\raux\cos\frac{2\pi p}{N}}{\rcyl^2+\raux^2-2\rcyl\raux\cos\frac{2\pi p}{N}}=B_p+1.
\end{equation}
Thus the right-hand side $D_p$ remains the same, while the new circulant matrix is obtained by adding $1$ to all matrix elements. It follows from \eqref{eq:DFT-relations} that $D'^{(m)}=D^{(m)}$ for all $m$; and from \eqref{eq:DFT-relations}, \eqref{sumexp}, and \eqref{Bm_exact} that 
\begin{equation}
\label{bmprime}
    B'^{(0)}=B^{(0)}+1=\frac{1}{1-\Big(\frac{\raux}{\rcyl}\Big)^N},\quad B'^{(m)}=B^{(m)},\quad m=1,2,\ldots,N-1. 
\end{equation}
These relations and \eqref{Im1} in turn imply
\begin{align}\label{Im-prime}
I'^{(m)}=I^{(m)}, \quad m=1,2,\ldots,N-1,
\end{align}
while, by \eqref{Im1}, \eqref{Dm_exact}, and \eqref{bmprime},
\begin{align}\label{I0-prime}
I'^{(0)}=\frac{D^{(0)}}{NB'^{(0)}}=\frac{I}{N}\frac{\Big(\frac{\rcyl}{\rfil}\Big)^{N}\bigg[1-\Big(\frac{\raux}{\rcyl}\Big)^N\bigg]}{1-\Big(\frac{\rcyl}{\rfil}\Big)^N}.
\end{align}
The MAS currents $I'_\ell$ are thus given by \eqref{Iell} or \eqref{Iell-cosines}, but with $I^{(0)}$ replaced by the $I'^{(0)}$ of \eqref{I0-prime}. Accordingly, \eqref{identity4} is modified to
\begin{align}
\label{identity4-prime}
    I'_\ell\sim\frac{I}{N}\left[\left(\frac{\rcyl}{\rfil}\right)^N+2\frac{(-1)^{\ell}t^{\frac{N+1}{2}}(t-1)\cos\frac{\pi\ell}{N}-t^2+t\cos\frac{2\pi\ell}{N}}{t^2-2t\cos\frac{2\pi\ell}{N}+1}\right], \quad \text{as}\,\,\, N\rightarrow\infty,
\end{align}
with $t$ as in \eqref{parameter-t}.
The conclusions of Sections \ref{asymptotic}--\ref{section-condition-number} are recast as follows.
\begin{itemize}
    \item{Case 1: Physical case $\raux>\rcri$}
    
    When $\raux>\rcri$, the asymptotic formula \eqref{asyILodd-tsmall} remains the same. The $I'_l$, when properly normalized, converge to (the same) surface current density $K_z(\phi)$, which can now be regarded as a (continuous) source located at $\rho=\raux$ and generating the scattered vector potential.  Since the only restriction on $\raux$ is $\rcri<\raux<\rcyl$, this source is certainly not unique. We can further show that $K_z(\phi)$ is the solution of an integral equation; we do not dwell on this, because the pertinent discussions within the context of electrodynamics \cite{fik_06} can easily be extended to the present case.

    \item{Case 2: Unphysical/oscillating case $\raux<\rcri$}
    
    In the more interesting case $\raux<\rcri$, we must neglect the exponentially small term $\left(\rcyl/\rfil\right)^N$ when simplifying \eqref{identity4-prime}. Accordingly, \eqref{asyILodd} is modified~to 
    \begin{equation}\label{asyILodd-prime}
    I'_{\ell}\sim\frac{2I}{N}\frac{(-1)^{\ell}t^{\frac{N+1}{2}}(t-1)\cos\frac{\pi\ell}{N}}{t^2-2t\cos\frac{2\pi\ell}{N}+1}, \quad \text{as}\,\,\, N\rightarrow\infty\quad (\raux<\rcri).
    \end{equation}
Consequently, the MAS currents $I'_\ell$ still exhibit exponentially large oscillations, but the mean is zero (previously, the mean was exponentially large). There is no associated continuous $K_z(\phi)$, so these $I'_\ell$ are unphysical in the same sense as were the $I_\ell$. 

\item{Case 3: Critical case $\raux=\rcri$}

Here, \eqref{identity4-prime} leads to $I'_0\sim I$ and, when $\ell\ne 0$, $I'_\ell\sim -I/N$. This time, the $\ell=0$ term in \eqref{Asc3} accounts for the image at the critical point $(\rcri,0)$---see the exact solution \eqref{Asc-exactsolution}---while, for $\ell\ne 0$, the normalized $I'_\ell$ amount to a $\phi$-independent surface-current density on the circle of radius $\raux=\rcri$, which replaces the image at the origin. 

\item In all three cases, it remains true that one obtains the true potential, i.e.,
\begin{equation}
\label{obtain-exact-solution-prime}
    A^{'\rm sc}_{z,N}\sim A_z^{\rm sc}, \quad \text{as}\,\,\, N\rightarrow\infty.
\end{equation}
This can be verified by going through the steps of \cref{convMAS1}. In particular, we have found (a posteriori) that the boundary condition at infinity is satisfied.

\item As $N\to\infty$, the condition number $\kappa'$ is large, but its order of magnitude differs. This is because the smallest $B^{(m)}$ is $B^{((N-1)/2)}$, while the largest is the $B'^{(0)}$ of \eqref{bmprime}---this is $O(1)$, while $B^{(0)}$ was exponentially small. Thus in place of \eqref{condition-number-kappa} we have
\begin{equation}
\label{condition-number-kappa-prime}
    \kappa'\sim\left(\frac{\rcyl}{\raux}\right)^\frac{N}{2}
\end{equation}
While exponentially large, this $\kappa'$ is much smaller than the $\kappa$ of \eqref{condition-number-kappa}.

\end{itemize}
As opposed to this second version of MAS, the first one had the advantage of a priori satisfying the boundary condition at infinity, for any finite $N$. Despite this, it presents the severe disadvantage of having a much larger matrix condition number. This affects the calculation of the auxiliary currents (when one finds them by solving the $N\times N$ system, something necessary for noncircular problems).  Furthermore, the currents of our first version oscillate about a mean value that is nonzero and large. This is disadvantageous, because stable algorithms should normally avoid large intermediate results \cite{higham}.

\section{Exterior Circular Problem: Numerical results}
\label{section:numerics}
\begin{figure}
\centering
\includegraphics[scale=0.51]{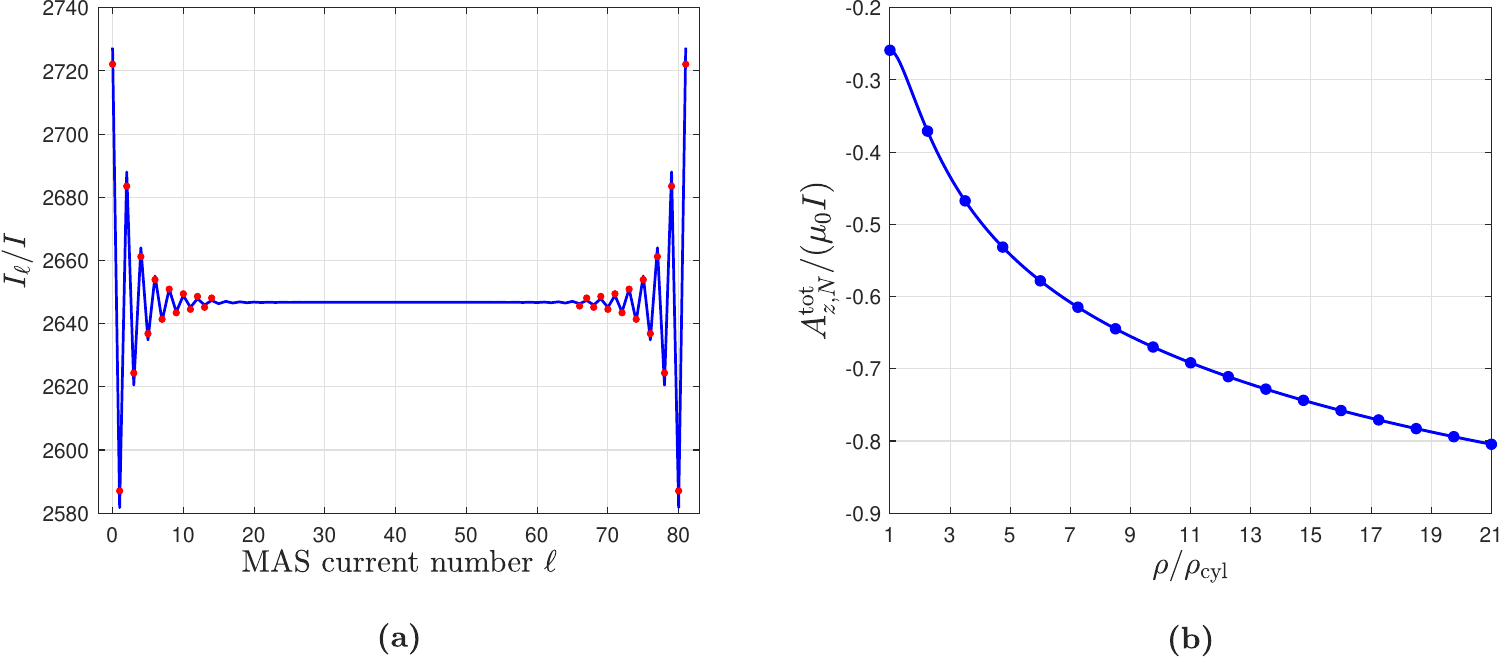}
\caption{Exterior circular problem with $\rcyl/\dref=8$, $\rfil/\dref=10$, $\raux/\dref=5.5$, and $N=81$. (a) Divergence of MAS currents: Normalized MAS currents $I_\ell/I$ vs the number $\ell$ of the auxiliary current. For symmetry reasons, the value of MAS current $I_0$ is repeated for $\ell=81$. Solid line: solution of the MAS system \eqref{sys1}; dots: asymptotic formula \eqref{asyILodd}.  (b) Convergence of MAS potential: Normalized total vector potential $A^{\rm tot}_{z,N}/(\mu_0 I)$ vs normalized observation distance $\rho/\rcyl$, for $\phi=45^{\circ}$. Solid line: MAS solution; dots: exact solution.}
\label{figure_numerical_1}
\end{figure}

Apart from some important numerical issues to be discussed in \cref{section:overdetermined}, \textit{all} theoretical predictions of the last two sections were verified numerically. In this section, we only illustrate some findings pertaining to Case~2, that is, when $\raux<\rcri$. \cref{figure_numerical_1} assumes $N=81$, $\rcyl/\dref=8$, and $\rfil/\dref=10$. The corresponding value of the critical radius is $\rcri/\dref=6.4$ and we take $\raux/\dref=5.5$. In the left figure, we have obtained a continuous curve by connecting the $I_\ell$-values by straight lines. The said values were obtained by employing bounded fundamental solutions and solving the $N\times N$ MAS system \cref{sys1} using the standard MATLAB solver. At least for the chosen parameters (more on this in \cref{section:overdetermined}), we obtain  identical results (at the scale of the figure) using the DFT solution \cref{Iell}--\cref{Iell-cosines}. The dots are corresponding results obtained from the asymptotic formula \cref{asyILodd}. There is very good agreement between the exact and asymptotic $I_\ell$, and we can clearly recognize the unphysical oscillations about the large mean value predicted asymptotically in \cref{asymptotic}. The right figure shows the total MAS potential obtained by substituting the oscillating $I_\ell$ into \cref{Asc1} (and adding the incident potential). This coincides, at the scale of the figure, with the exact total potential \cref{Asc-exactsolution}. Thus MAS gives the correct potential despite the unphysical oscillations.

For brevity, we show no results pertaining to traditional fundamental solutions. But, as predicted theoretically, there arises only one essential difference with \cref{figure_numerical_1}, namely the mean value of the oscillations, which is now close to zero.

\section{Extensions}
\label{section-extensions}

We have discussed several aspects of the application of MAS to exterior circular Laplace-Neumann problems. We now extend our results in a number of different directions. 

\subsection{Purely numerical effects; 
 the use of overdetermined systems}
\label{section:overdetermined}

Errors due to roundoff can in principle be alleviated by longer computer wordlengths. As we already discussed, the unphysical oscillations analyzed herein have nothing to do with roundoff; the  existence of an asympotic formula predicting the oscillating values (see the left \cref{figure_numerical_1}) suffices to confirm this observation, which is also supported by the agreement---already mentioned in \cref{section:numerics}---between the two methods of obtaining the $I_\ell$ (namely, DFT and direct solution of the $N\times N$ system). For the parameters in \cref{figure_numerical_1}, this agreement is not implausible, because the condition number is not prohibitively large: \cref{condition-number-kappa} gives $\kappa=5\cdot10^{12}$ (the MATLAB routine \textit{cond} gives the same $\kappa$). 

While modern computers can deal with systems having such condition numbers, $\kappa$ rapidly increases with $N$. For example, if we increase the $N$ of \cref{figure_numerical_1} to $N=101$, we have $\kappa\sim 10^{16}$, and the two methods of obtaining $I_\ell$ no longer give the same results: DFT gives $I_0/I=4\cdot10^4$, while directly solving \cref{sys1} gives $I_0/I=1\cdot10^4$. Now, roundoff (which was negligible in \cref{figure_numerical_1}) becomes important (at least in the calculation of $I_\ell$, recall the insensitivity of \cref{section-insensitivity}). 

The MAS-related literature (see for example \cite{che_hon_20,eremin-tsitsas-kouroublakis-fikioris}) often uses overdetermined systems. In our case, such schemes amount to enforcing
the boundary condition \cref{bc0} at $M$ collocation points $(\rcyl,2\pi p/M)$ ($p=0,\dots,M-1$ and $M>N$) and solving the resulting overdetermined $M\times N$ system (for $I_\ell$) in the least-squares sense. We experimented with such schemes, and found no significant improvement as far as the oscillations in the unphysical case are concerned. 

Let us return to square ($N\times N$) systems. For sufficiently large $N$ (or sufficiently small $\raux$), we found that the  $I_\ell$ obtained by solving the system no longer agree with the asymptotic formula \cref{asyILodd}. This means that solution behavior is clouded by roundoff: the \textit{true} behavior (which would be obtained in a hypothetical computer with infinite wordlength) is, somewhat paradoxically, better predicted by the asymptotic formula. In other words, the effects caused by the large condition numbers are purely numerical; so are the mitigating effects due to insensitivity. These two effects depend on computer wordlength and the system solver and, more generally, on the specific software and hardware one uses.

\subsection{A noncircular problem}\label{noncircular_external_problem}
In practice, MAS is applied to noncircular problems \cite{tsi_zou_fik_lev_18,dou_zha_li_che_20}. To demonstrate that our main findings carry over to noncircular shapes, we consider an elliptic cylinder of infinite permeability (as in \cref{geometry}, but with the circular cylinder replaced by an elliptic one). This problem presents the advantage of having a particularly simple exact solution; it is found in \cref{appendix-elliptic cylinder exact solution} using separation of variables in conjunction with the integral $J$ in \cref{idefinition} and \cref{ievaluated}.

The ellipse's major (minor) axis lies along the $x$- ($y$-) axis and has length $2a$ ($2b$). The current filament $I$ is again located at $(\rho,\phi)=(\rfil,0)$ with $\rfil>a$. The Neumann boundary condition on the surface of the ellipse is 
\begin{align}\label{Neumann_bc_elliptic}
\frac{\partial A^{\rm sc}_z}{{\partial n}}=-\frac{\partial A^{\rm inc}_z}{{\partial n}},
\end{align}
where $\partial/{\partial n}$ is the normal derivative; an expression for the latter can be readily obtained by using the polar equation of the ellipse $\rho=b/(1-h^2\cos^2\phi)^{1/2}$, $-\pi<\phi\leq \pi$, where $h=(1-b^2/a^2)^{1/2}$ is the eccentricity.

We  choose the MAS currents $I_\ell$ to lie on a confocal elliptic auxiliary surface with major and minor axes $2a_{\rm aux}$ and $2b_{\rm aux}$, respectively ($a_{\rm aux}<a$ and $b_{\rm aux}<b$). The confocal property allows for both the physical and the auxiliary surface to be readily described in the same elliptic coordinate system. In this context, the currents $I_\ell$ are uniformly distributed upon the auxiliary surface with respect to the elliptic angular coordinate $\eta$, i.e., at $\eta_\ell=2\pi\ell/N$ $(\ell=0,\dots,N-1)$. Then, the polar coordinates $(\rho_\ell,\phi_\ell)$ of $I_\ell$ can be straightforwardly obtained and, 
depending on the fundamental solution we employ, the MAS scattered potential is given by either \cref{Asc1} or \cref{Asc3}, with  
\begin{align}\label{Rl_elliptic} 
R_\ell=\sqrt{\rho^2+\rho_\ell^2-2\rho\rho_\ell\cos(\phi-\phi_\ell)}.  
\end{align}

\begin{figure}
\centering
\includegraphics[scale=0.51]{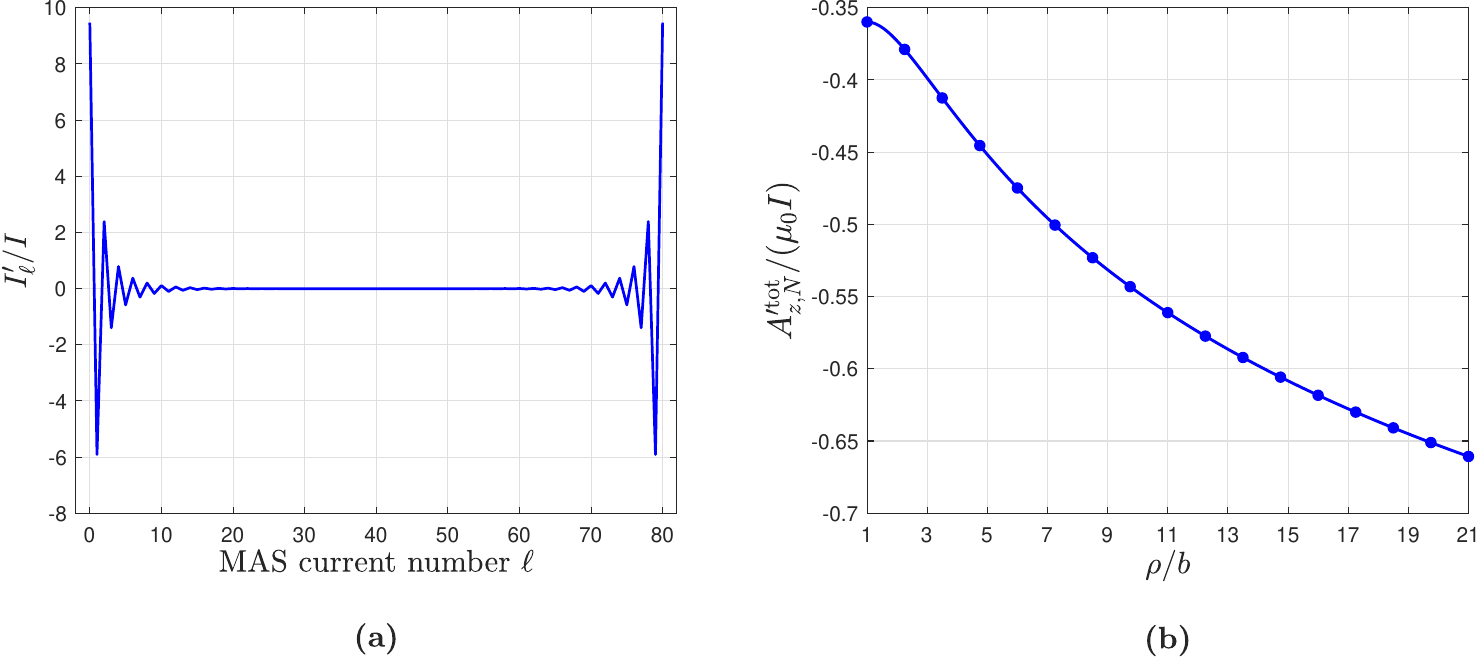}
\caption{Exterior elliptic problem with $a/\dref=6$, $b/\dref=3$, $\rfil/\dref=7.5$, $a_{\rm aux}/\dref=5.2222
$, $b_{\rm aux}/\dref=0.5205$, and $N=80$. (a) Divergence of MAS currents: Normalized MAS currents $I'_\ell/I$ vs the number $\ell$ of the auxiliary current. For symmetry reasons, the value of MAS current $I'_0$ is repeated for $\ell=80$.  (b) Convergence of MAS potential: Normalized total vector potential $A'^{\rm tot}_{z,N}/(\mu_0 I)$ vs normalized observation distance $\rho/b$, for $\phi=90^{\circ}$. Solid line: MAS solution; dots: exact solution.}
\label{figure_elliptic}
\end{figure}

Similarly, the $N$ collocation points are uniformly placed upon the surface of the elliptic cylinder with respect to $\eta$, i.e., at $\eta_p=2\pi p/N$ $(p=0,\dots,N-1)$, and their polar coordinates $(\rho_p,\phi_p)$ can be obtained. Then, \eqref{Neumann_bc_elliptic} is enforced at $(\rho_p,\phi_p)$ and we finally arrive at a $N\times N$ system of linear equations for the determination of $I_\ell$. Unlike the circular case, now the system is not circulant and cannot be solved with DFTs.

We show numerical results in \cref{figure_elliptic}, where the line source $I$ is located at $\rfil/\dref=7.5$ and illuminates an elliptic cylinder with $a/\dref=6$ and $b/\dref=3$. The auxiliary surface is a confocal ellipse with $a_{\rm aux}/\dref=5.2222$ and $b_{\rm aux}/\dref=0.5205$. MAS is employed with the logarithmic fundamental solutions of \cref{MAStrad} (and therefore the notation $A'_z$ and $I'_\ell$ is used instead of $A_z$ and $I_\ell$). The left figure shows $I'_\ell/I$. We see that the currents diverge, exhibiting the familiar oscillation pattern with maximum oscillations observed near $\phi=0$. The right figure shows the plot of the total MAS potential. The values corresponding to the solid line are obtained by means of \cref{Az_inc}, \cref{Asc3}, and \cref{Rl_elliptic}, using the oscillating $I'_\ell$ of the left figure. 

Values of the aforementioned exact solution (of \cref{appendix-elliptic cylinder exact solution}) are also shown with dots. At the scale of the figure, there is perfect agreement despite the oscillations. Similar observations hold if MAS is employed with the bounded fundamental solutions $\ln(R_\ell/\rho)$: Again, diverging MAS currents yield the correct vector potential, although the oscillations of $I_\ell$ now exhibit more unusual patterns.

Needless to say, whenever oscillations were observed, the auxiliary surface did not enclose the singularities of the scattered potential. These singularities can be easily determined from the (exterior) solution found in \cref{appendix-elliptic cylinder exact solution}, which has an obvious analytic continuation to the interior.

\subsection{Interior circular problem (cavity problem)}
\label{internal_problem}
\begin{figure}
\centering
\includegraphics[scale=0.47]{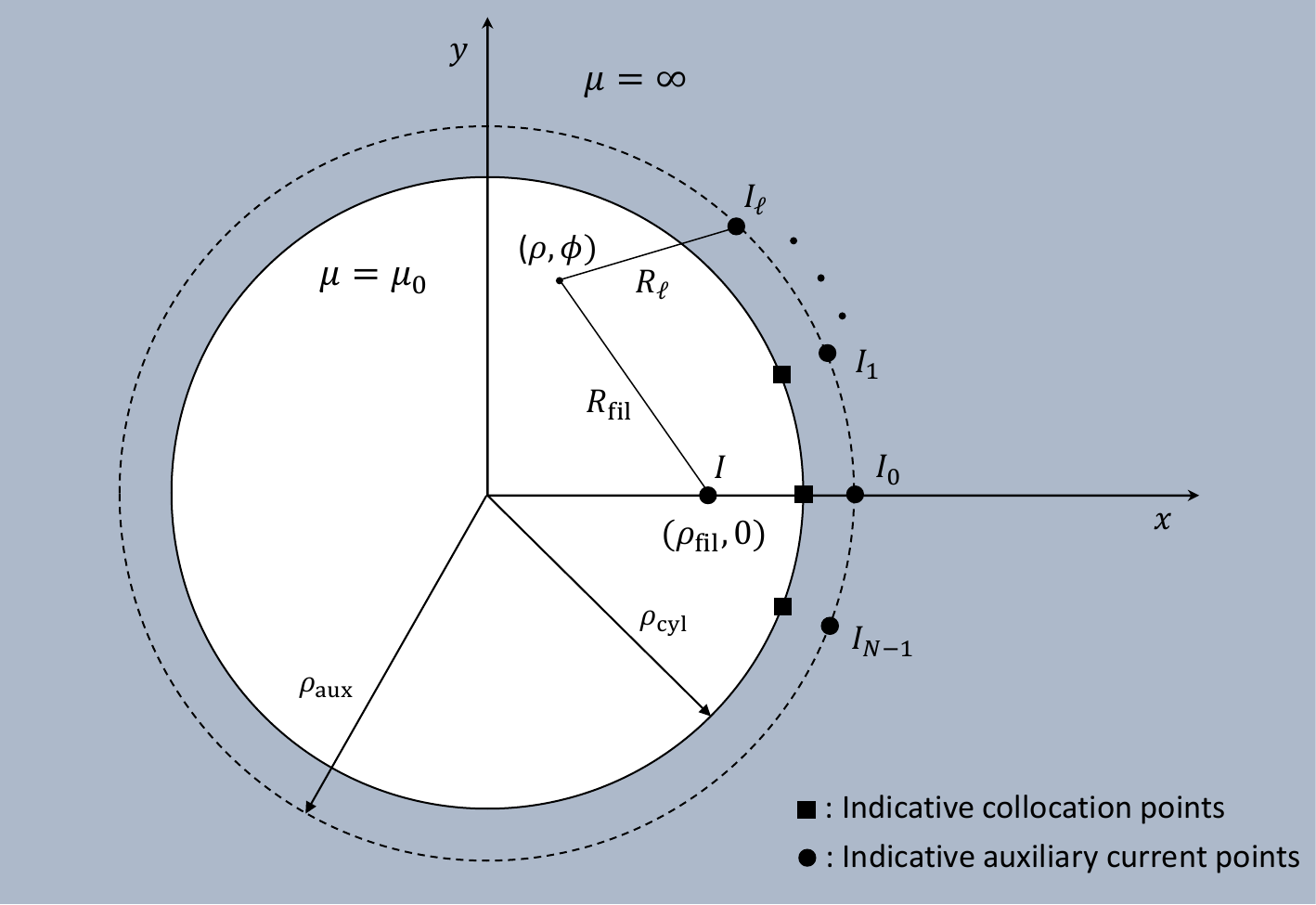}
\caption{Geometry of the interior circular problem. Indicative positions of MAS currents and collocation points are depicted.}
\label{internal_geometry}
\end{figure}

Now consider the circular ``cavity'' problem of \cref{internal_geometry}, in which a line current $I$ is located at $(\rho,\phi)=(\rfil,0)$ \textit{inside} a nonmagnetic ($\mu=\mu_0$) cylinder of radius $\rcyl$, with the surrounding medium having infinite permeability ($\mu=\infty$).  As in \cref{Problemdesciption}, no current can flow at $\rho>\rcyl$. However, the properties of $\mathbf{H}$ and $\mathbf{B}$ now differ. In the surrounding medium, with $\mu=\infty$,  $\mathbf{H}$ must be finite and so $\mathbf{B}=\mu\mathbf{H}$ becomes infinite. In order for the tangential components of $\mathbf{H}$ and the normal components of $\mathbf{B}$ to be continuous across the surface of the cylinder $\rho=\rcyl$ (boundary conditions), the field lines of $\mathbf{H}$ (and $\mathbf{B}$) should be tangential at $\rho=\rcyl+0$, but not at $\rho=\rcyl-0$. In the interior, where $\mu=\mu_0$, the direction of the field lines is not normal to the surface of the cylinder (as in the exterior circular problem), but is such that the aforementioned boundary conditions be valid there \cite{haus_melcher,bladel}. The component $H_{\phi}$ of $\mathbf{H}$  equals $I/(2\pi\rho)$ for $\rho\geq\rcyl$ \cite{bladel}. Therefore the total interior potential satisfies  $\partial A_{z}^{\rm tot}/\partial\rho=-\mu_0I/(2\pi\rcyl)$ at $\rho=\rcyl-0$.   

Once again, we formulate our problem in terms of the scattered potential $A_z^\mathrm{sc}$. Thus $A_z^\mathrm{sc}$ satisfies Laplace's equation in the interior, together with the Neumann boundary condition
\begin{equation}
\label{neumann-internal-circular}
    \frac{\partial A_{z}^{\rm sc}}{\partial\rho}=-\frac{\partial A_{z}^{\rm inc}}{\partial\rho}-\frac{\mu_0I}{2\pi\rcyl},\,\,\, \text{at}\,\,\ \rho=\rcyl-0,
\end{equation}
where $A_{z}^{\rm inc}$ is given by \eqref{Az_inc} and \eqref{Rfil} (but now, $\rfil<\rcyl$). 
By Amper\'e's circuital law, 
$\oint_{\rho=\rcyl-0}
\mathbf{H}^\mathrm{inc}\cdot{\rm d}\mathbf{l}=I$, so that 
\begin{equation}
\label{self-consistency-internal-circular}
\int_0^{2\pi}\left(-\frac{\partial A_{z}^{\rm inc}}{\partial\rho}-\frac{\mu_0I}{2\pi\rcyl}\right)\rcyl\,d\phi=0,\,\,\, \text{at}\,\,\ \rho=\rcyl-0,
\end{equation}
meaning that the Neumann boundary data (i.e., the right-hand side of \eqref{neumann-internal-circular}) satisfies the self-consistency condition.  For an independent verification of \eqref{self-consistency-internal-circular}, evaluate its left-hand-side integral using \eqref{Az_inc}, \eqref{Rfil}, and the bottom integral \eqref{ievaluated-circle-m-equals-zero}.  

The (non-unique) solution to our boundary-value problem can be found using elementary methods to be
\begin{align}\label{Ascinterior}
A_{z}^{\rm sc}(\rho,\phi)=A_0-\frac{\mu_0I}{2\pi}\ln\frac{\sqrt{\rho^2+\rcri^2-2\rho\rcri\cos\phi}}{\rcri}, \quad \rho\le \rcyl,
\end{align}
where $A_0$ is an arbitrary constant, i.e., $A_0$ is independent of the observation coordinates $\rho$ and $\phi$. In \eqref{Ascinterior}, $\rcri$ is defined as before, but note that the various distances now satisfy
\begin{equation}
    \rfil<\rcyl<\raux,\quad \rcyl<\rcri=\frac{\rcyl^2}{\rfil}.
\end{equation}

Our purpose is to extend our previous analysis (of the exterior problem) and, especially, our previous investigation of the unphysical case. To this end, we apply a MAS scheme that is a straightforward analogue of the method in \cref{MAStrad} (i.e., use traditional fundamental solutions; there seems to be no  analogue to the method of \cref{MASkats}, because the $\ln(R_\ell/\rho)$ are  singular at $\rho=0$). We thus place $N$ auxiliary sources $I_\ell$ at $(\rho,\phi)=(\raux, 2\pi\ell/N)$ ($\ell=0,1,\ldots,N$) write 
\begin{align}\label{AscNin}
A_{z,N}^{\rm sc}(\rho,\phi)=-\frac{\mu_0}{2\pi}\sum_{\ell=0}^{N-1}I_\ell\ln\frac{\sqrt{\rho^2+\raux^2-2\rho\raux\cos\Big(\phi-\frac{2\pi\ell}{N}\Big)}}{\dref},\quad \rho\le\rcyl,
\end{align}
and determine  
the $I_\ell$ by enforcing \cref{neumann-internal-circular} at $(\rho,\phi)=(\rcyl, 2\pi p/N)$ ($p=0,1,\ldots,N$). The resulting circulant $N\times N$ system, in which $\dref$ and $A_0$ do not appear, can be solved just as in \cref{MAStrad}. The final result is given by \eqref{Iell} or \eqref{Iell-cosines}, where now 
\begin{align}
I^{(0)}&=\frac{I}{N}\frac{\Big(\frac{\rfil}{\rcyl}\Big)^{N}\bigg[1-\Big(\frac{\rcyl}{\raux}\Big)^N\bigg]}{\Big(\frac{\rcyl}{\raux}\Big)^{N}\bigg[1-\Big(\frac{\rfil}{\rcyl}\Big)^N\bigg]},\label{I0_internal_exact} \\
I^{(m)}&=\frac{I}{N}\frac{\bigg[\Big(\frac{\rfil}{\rcyl}\Big)^{N-m}+\Big(\frac{\rfil}{\rcyl}\Big)^m\bigg]\bigg[1-\Big(\frac{\rcyl}{\raux}\Big)^N\bigg]}{\bigg[\Big(\frac{\rcyl}{\raux}\Big)^{N-m}+\Big(\frac{\rcyl}{\raux}\Big)^m\bigg]\bigg[1-\Big(\frac{\rfil}{\rcyl}\Big)^N\bigg]}, \quad m=1,\ldots,N-1. \label{Im_internal_exact}
\end{align}
In \eqref{AscNin}, apply the identity \eqref{identity3} together with the definition \eqref{eq:DFT-relations} of the DFT to obtain the exact relation
\begin{align}\label{AscN1inexact}
A_{z,N}^{\rm sc}(\rho,\phi)=-\frac{\mu_0}{2\pi}NI^{(0)}\ln\frac{\raux}{\dref}+\frac{\mu_0}{4\pi}\sum_{m\neq 0}\frac{NI^{(m)}}{|m|}\bigg(\frac{\rho}{\raux}\bigg)^{|m|}e^{im\phi},\quad \rho\le\rcyl,
\end{align}
in which the $I^{(m)}$ are the $N$-periodic and even (in $m$) extensions of the $I^{(m)}$ in \eqref{Im_internal_exact} (this renders the series in \eqref{AscN1inexact} convergent). The following large-$N$ formulas are simple consequences of \cref{I0_internal_exact} and \cref{Im_internal_exact}, 
\begin{align}
I^{(0)}&\sim\frac{I}{N}\bigg(\frac{\raux}{\rcri}\bigg)^N,\quad m=0,\quad \mathrm{as\ } N\to\infty, \label{asyI0_internal}\\
I^{(m)}&\sim
\frac{I}{N}\bigg(\frac{\raux}{\rcri}\bigg)^{|m|},\quad m=\pm 1,\pm 2,\ldots,\pm\frac{N-1}{2}, 
\quad \mathrm{as\ } N\to\infty.
\label{asyIm_internal}
\end{align}
In \cref{asyIm_internal}, we used $|m|$ in the exponent so as to extend  to negative $m$.
We obtain a large-$N$ formula for the \textit{second term} in  \eqref{AscN1inexact} (i.e., the series) by replacing $I^{(m)}$ with the asymptotic expression in \cref{asyIm_internal}; the resulting series, which is independent of $\raux$, can be summed with the aid of \cref{identity3}:
\begin{multline}\label{AscN1insecondterm}
\frac{\mu_0}{4\pi}\sum_{m\neq 0}\frac{NI^{(m)}}{|m|}\bigg(\frac{\rho}{\raux}\bigg)^{|m|}e^{im\phi}\sim
\frac{\mu_0 I}{4\pi}\sum_{m\neq 0}\frac{1}{|m|}\bigg(\frac{\rho}{\rcri}\bigg)^{|m|}e^{im\phi} \\ =-\frac{\mu_0I}{2\pi}\ln\frac{\sqrt{\rho^2+\rcri^2-2\rho\rcri\cos\phi}}{\rcri}, \quad \mathrm{as\ } N\to\infty.
\end{multline}
\begin{figure}
\centering
\includegraphics[scale=0.72]{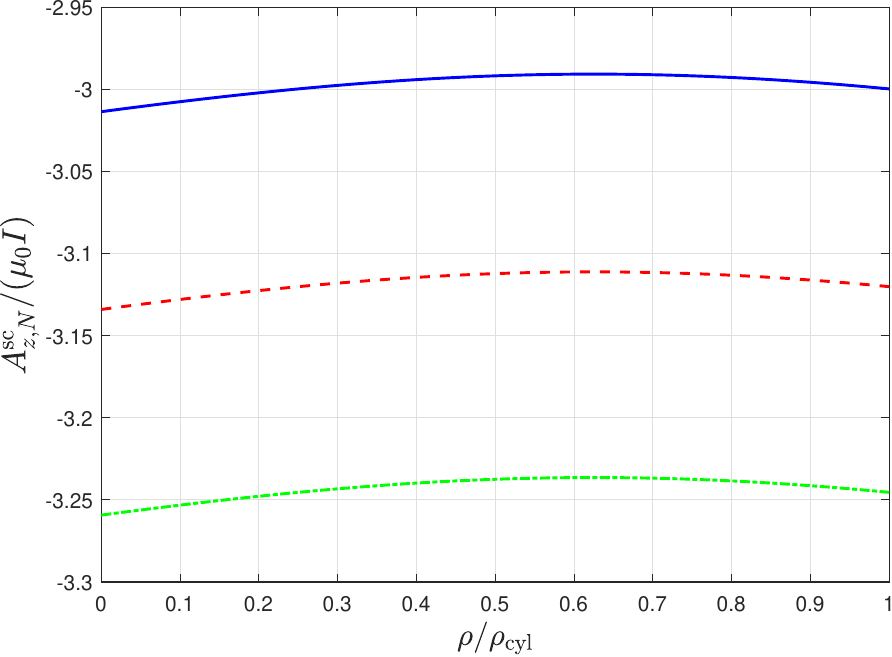}
\caption{Interior circular problem with $\rcyl/\dref=5$, $\rfil/\dref=4$, and $\raux/\dref=6.5$. Divergence of MAS potential: Normalized scattered vector potential $A^{\rm sc}_{z,N}/(\mu_0 I)$ (as calculated from \eqref{AscNin}) vs normalized observation distance $\rho/\rcyl$, for $\phi=60^{\circ}$. Blue solid line: $N=59$; red dashed line: $N=60$; green dash-dotted line: $N=61$.}
\label{figure_numerical_2}
\end{figure}
By \cref{asyI0_internal}, the \textit{first term} in \cref{AscN1inexact} is exponentially small when $\raux<\rcri$, but exponentially large when $\raux>\rcri$. We are led to the following conclusions.

\begin{itemize}
    \item {Case 1: Physical case $\raux<\rcri$} 
    
    When $\raux<\rcri$, the normalized $I_\ell$ converge to a surface current density $K_z(\phi)$, just as in the interior problem. Furthermore (because the first term in \cref{AscN1inexact} converges to zero as $N\to\infty$), $A_{z,N}^{\rm sc}$ converges to a true solution of the problem, specifically to the right-hand side of \cref{Ascinterior} with $A_0=0$.
    
    \item {Case 2: Unphysical/oscillating case $\raux>\rcri$} 
    
    Here, the $I_\ell$ cannot converge to a surface current density; rather, their asymptotic behavior is given by \cref{asyILodd}, but with $t=\raux/\rcri>1$. Furthermore---with an exception to be discussed below---$A_{z,N}^{\rm sc}$ diverges as $N\to\infty$. The ``divergent part''---i.e., the first term in  \eqref{AscN1inexact}---is exponentially large in $N$, and independent of $\rho$ and $\phi$, but is dependent on $\raux$. 
\end{itemize}

For Case 2, it is worth stressing that the mean value of the oscillating $I_\ell$ is large, as in the exterior problem with the \textit{bounded} fundamental solutions, even if we  used the \textit{traditional} ones here. Yet what is really new is the divergence of the potential, something not encountered in any of the two schemes we applied to the exterior problem. We illustrate this divergence in \cref{figure_numerical_2}, where it is seen that three successive values of $N$ give potentials that (due to the change of the first term in \cref{AscN1inexact}) differ by a $\rho$-independent constant. 

In Case~2, choosing $\raux=\dref$  has the desirable consequence of rendering the divergent part to be zero; this is the exceptional case we previously alluded to. It seems, however, that this case is only of theoretical interest because there is no straightforward way of extending the choice $\raux=\dref$ to noncircular problems.  

It is perhaps more interesting that the divergent part is independent of the observation point $(\rho,\phi)$:  If it is $\mathbf{H}$ that is finally desired, and if we calculate $\mathbf{H}$ numerically from $A_{z,N}^{\rm sc}$, then all sufficiently large values of $N$ will give the same $\mathbf{H}$. (As large intermediate results should in general be avoided \cite{higham}, this is true with the proviso that $A_{z,N}^{\rm sc}$ does not become too large.) Thus in Case 2, the first term in \eqref{AscN1inexact} somewhat resembles a nonzero $A_0$ in the exact solution \eqref{Ascinterior}; but (since it is exponentially large) the said term does not \textit{converge} to an $A_0$. 

These findings render the MAS behavior of the present (interior circular Laplace-Neumann) problem distinctly different from that of the electrodynamics problems we have studied \cite{fik_06, fik_psar_07, fik_tsi_15, tsi_zou_fik_lev_18, fikioris2018bookchapter}.

\section{Conclusions and future work}
\label{sec:conclusions}

Previous works \cite{fik_06, fik_psar_07, fik_tsi_15, tsi_zou_fik_lev_18, fikioris2018bookchapter} have shown that it is possible for the final MAS solution to converge to the correct one despite the divergence and unphysical oscillations of an intermediate result, namely the MAS currents. Here, we extended these findings to the case of Laplace-Neumann problems. As in the previous works, the normalized MAS currents diverge and oscillate whenever the analytic continuation of the scattered potential presents a singularity between the auxiliary surface and the 2D boundary; in the circular problems, the oscillations are well-predicted by asymptotic formulas (\eqref{asyILodd} and \eqref{asyILodd-prime}) and would therefore occur even with perfect hardware and software; these oscillations must be distinguished from effects (\cref{section:overdetermined}) that are purely numerical in nature and depend on the hardware/software implementation.

The present analysis (for the Laplace-Neumann case) has certain particularities that do not pertain to the corresponding time-harmonic cases: The exact and asymptotic formulas are generally simpler, so that conclusions arise in a clearer way. We treated the exterior circular problem by two different methods and found noteworthy differences between the two arising unphysical/oscillating cases (in one of the two, the currents oscillate about a \textit{large} mean value, whereas the mean is zero in the corresponding time-harmonic problems). The interior circular magnetostatic problem, which is known to require a different Neumann formulation from the exterior one (compare \cref{bc0} to \cref{neumann-internal-circular}), gave rise to a ``divergent'' vector potential (in the sense explained in \cref{internal_problem}; no such divergence seems to occur in the corresponding time-harmonic problems). Finally, to extend our findings to noncircular problems, we chose an elliptic boundary, for which we developed a particularly simple closed-form solution (\cref{appendix-elliptic cylinder exact solution}), whose singularities can be deduced a priori.  

The ``insensitivity'' of \cref{section-insensitivity} was observed numerically. In fact, the numerical results suggest that the present (Laplace-Neumann) problem exhibits more insensitivity than certain time-harmonic problems we have studied \cite{fik_06, fik_psar_07, fik_tsi_15, tsi_zou_fik_lev_18, fikioris2018bookchapter}. For now, we make no attempt to state this empirical observation in a more precise manner, but we do plan to study it in the future. Such a study  should consider insensitivity concurrently with the often overlooked phenomenon of oscillations (which, as we saw in \cref{section-insensitivity}, can be thought of as a competing effect).
We also plan to study more carefully the manner in which the two combined phenomena affect the MFS solutions of noncircular problems.

\appendix
\section{Evaluation of the integral $J(m,x,y)$}
\label{appendix-integral}
In this appendix we start from \cref{idefinition} and prove \cref{ievaluated}. It is convenient to define 
\begin{align}
\label{parameter-tau}
\tau=\sqrt{x^2+y^2-1}.
\end{align}
The inequalities in \cref{idefinition} give
\begin{align}
\label{parameter-tau-inequality}
0<\tau\le x. 
\end{align}
The definitions \cref{idefinition} and \cref{parameter-tau} yield
\begin{align}
\label{ievaluation-1}
J=\int_{-\pi}^{\pi}\frac{(1-x\cos\theta-i\tau\sin\theta)\cos(m\theta)}{(1-x\cos\theta)^2+\tau^2\sin^2\theta}{\rm d}\theta.
\end{align}
(The additional term $-i\tau\sin\theta$, which is odd in $\theta$, does not contribute to the integral.) Therefore
\begin{align}
\label{ievaluation-2}
J=\int_{-\pi}^{\pi}\frac{\cos(m\theta)}{1-x\cos\theta+i\tau\sin\theta}{\rm d}\theta.
\end{align}
Set $z=e^{i\theta}$ to express $J$ as a  contour integral,
\begin{align}
\label{ievaluation-3}
J=i\oint_{|z|=1}\frac{z^{2m}+1}{z^m[(x-\tau)z^2-2z+(x+\tau)]}{\rm d}z.
\end{align}
As long as $y\ne 1$, \cref{parameter-tau} implies that $x\ne \tau$ and that the discriminant of the quadratic in \cref{ievaluation-3} equals $4y^2$. Thus
\begin{align}
\label{ievaluation-4}
J=\frac{i}{x-\tau}\left[\oint_{|z|=1}\frac{z^{2m}}{z^m(z-z_1)(z-z_2)}{\rm d}z+\oint_{|z|=1}\frac{{\rm d}z}{z^m(z-z_1)(z-z_2)}\right],
\end{align}
where the two roots $z_1$ and $z_2$ are  given by
\begin{align}  z_{1,2}=\frac{1\mp y}{x-\tau}=\frac{x+\tau}{1\pm y},
\end{align}
in which the last expression  followed from \cref{parameter-tau}. By \cref{parameter-tau-inequality}, the inequalities in \cref{idefinition}, and the assumption $y\ne 1$, the roots are positive, unequal, and lie outside the unit circle. Consequently the first integral in \cref{ievaluation-4} vanishes, whereas the integrand of the second has one singularity within the unit circle, namely a pole of order $m$ at $z=0$. The residue there can be found by means of the series
\begin{align}
\label{ievaluation-10}
\frac{1}{z^m(z-z_1)(z-z_2)}=\frac{1}{z_2-z_1}
\sum_{n=0}^\infty\left(z_1^{-n-1}-z_2^{-n-1}\right)z^{n-m}.
\end{align}
Application of the residue formula and re-introduction of $\tau$ from \cref{parameter-tau} then yields \cref{ievaluated}. By analytic continuation, \cref{ievaluated} is also valid for the excluded case $y=1$.

\section{Circular problems: Solution of circulant systems}
\label{appendix-circulant}
In this appendix, we show \cref{Bm_exact}--\cref{Im_exact}. As discussed in our Introduction, we first apply the DFT to the symmetric circulant system \eqref{sys1}. We obtain 
\begin{align}\label{Im1}
    I^{(m)}=\frac{D^{(m)}}{NB^{(m)}}.
\end{align}
where $B^{(m)}=B^{(N-m)}$ is the DFT of $B_p$ (i.e., the first row of the system-matrix of \cref{sys1}), $D^{(m)}=D^{(N-m)}$ is the DFT of the right-hand side vector of \cref{sys1}, and $I^{(m)}=I^{(N-m)}$ is the DFT of the solution $I_\ell$, where the DFTs are defined according to \cref{eq:DFT-relations}. 
We now substitute \eqref{Dp} into the equation defining $D^{(m)}$ and apply the identity \cref{identity1} with $\rho_1=\rcyl$ and $\rho_2=\rfil$ (see \cref{distance-inequality}) to obtain
\begin{align}\label{Dm1}
 D^{(m)}&=\frac{1}{N}\sum_{p=0}^{N-1}\Bigg[-I\frac{\rcyl^2-\rcyl\rfil\cos\frac{2\pi p}{N}}{\rcyl^2+\rfil^2-2\rcyl\rfil\cos\frac{2\pi p}{N}}\Bigg]e^{-i2\pi mp/N} \notag\\
 &\stackrel{}{=}\frac{-I\rcyl}{N}\sum_{p=0}^{N-1}\Bigg[-\frac{1}{2\rcyl}\sum_{n\neq 0}\bigg(\frac{\rcyl}{\rfil}\bigg)^{|n|}e^{i2\pi np/N}\Bigg]e^{-i2\pi mp/N}  \notag\\
 &=\frac{I}{2}\sum_{q=-\infty}^{\infty}\bigg(\frac{\rcyl}{\rfil}\bigg)^{|qN+m|},\quad \text{with}\,\,q\neq0\,\,\text{if}\,\,m=0,
\end{align}
in which the last step followed by interchanging the order of summation and applying the identity \cref{sumexp}. With the aid of the geometric series we then obtain
\begin{align}\label{Dm_exact}
D^{(0)}=I\frac{\Big(\frac{\rcyl}{\rfil}\Big)^N}{1-\Big(\frac{\rcyl}{\rfil}\Big)^N}, \quad D^{(m)}=\frac{I}{2}\frac{\Big(\frac{\rcyl}{\rfil}\Big)^{N-m}+\Big(\frac{\rcyl}{\rfil}\Big)^m}{1-\Big(\frac{\rcyl}{\rfil}\Big)^N},\quad m=1,\dots,N-1, 
\end{align}
Similarly, starting from \cref{eq:DFT-relations} and \cref{Bl}, we use \eqref{identity2} (with $\rho_1=\rcyl$ and $\rho_2=\raux$---see \cref{distance-inequality}) and \cref{sumexp} to obtain
\begin{align}\label{Bm1}
 B^{(m)}=\frac{1}{2}\sum_{q=-\infty}^{\infty}\bigg(\frac{\raux}{\rcyl}\bigg)^{|qN+m|},\quad \text{with}\,\,q\neq0\,\,\text{if}\,\,m=0,
\end{align}
from which \cref{Bm_exact} follows. Equations \cref{Iell}--\cref{Im_exact} then follow from \cref{eq:DFT-relations}, \cref{Im1}, \cref{Dm_exact}, and \cref{Bm_exact}.

\section{Exact solution for elliptic cylinder}\label{appendix-elliptic cylinder exact solution}
Here, we outline the derivation of the exact solution to the 2D problem of an elliptic cylinder of infinite permeability, illuminated by a line-current source $I$. The cross section of the cylinder is an ellipse with major and minor axes $2a$ and $2b$, respectively. The current $I$ is located at $(x,y)=(\rfil,0)$. The problem is solved subject to the Neumann boundary condition for the total vector potential at the surface of the cylinder. 
We employ the separation of variables in elliptic coordinates to obtain the solution in a particularly simple form. 
\par
The 2D elliptic coordinates $(\xi,\eta)$ are connected to the Cartesian coordinates $(x,y)$ by means of the transformation
\begin{align}\label{elliptic_coord}
x=c\cosh\xi\cos\eta, \quad y=c\sinh\xi\sin\eta,
\end{align}
where $\xi\geq0$, $-\pi<\eta\leq\pi$, and $c$ is a parameter that determines the focal points of the elliptic coordinate system. In the latter, the surfaces of constant $\xi$ are confocal ellipses, while  the surfaces of constant $\eta$ are branches of confocal hyperbolas; in both cases the focal distance is $2c$. 

In elliptic coordinates, the surface of the elliptic cylinder is described by $\xi=\xi_0$, where $\xi_0$ is determined by either $\cosh\xi_0=a/c$ or $\sinh\xi_0=b/c$, with $c=(a^2-b^2)^{1/2}$ half the focal distance. The respective coordinates of the current filament are $(\xi,\eta)=(\xifil,0)$, where $\cosh{\xifil}=\rfil/c$. For $\xi>\xi_0$ and $(\xi,\eta)\neq(\xifil,0)$, $A_z^{\rm sc}$ satisfies the 2D Laplace equation. In accordance with \eqref{bc0}, the Neumann boundary condition is 
\begin{align}
\label{neumann-external-elliptic}
\frac{1}{h_\xi}\frac{\partial A_z^{\rm sc}}{\partial\xi}=-\frac{1}{h_\xi}\frac{\partial A_z^{\rm inc}}{\partial\xi},\,\,\, \text{at}\,\,\,\xi=\xi_0+0.
\end{align}
where $h_\xi=c(\cosh^2\xi-\cos^2\eta)^{1/2}$ is the appropriate scale factor.
Expression \eqref{Az_inc} for $A_z^{\rm inc}$ is written in elliptic coordinates as
\begin{align}
\label{Azinc_elliptic_coords}
    A_z^{\rm inc}(\xi,\eta)=-\frac{\mu_0I}{2\pi}\ln\frac{\sqrt{(c\cosh\xi\cos\eta-\rfil)^2 + c^2\sinh^2\xi\sin^2\eta}}{\dref},
\end{align} 
and we can easily verify that the Neumann boundary data (i.e., the right-hand side of \eqref{neumann-external-elliptic}) satisfies the self-consistency condition: Indeed, the related integral is reduced to \eqref{idefinition} with $m=0$, $x=\rfil/a$, and $y=b/a$, and through \eqref{ievaluated} turns out to vanish.

Because the Laplace equation is separable in elliptic coordinates, the solution for $A_z^{\rm sc}$ can be written in the form of the series
\begin{align}
\label{seriesforaz}
    A_z^{\rm sc}(\xi,\eta)=-\frac{\mu_0I}{2\pi}\sum_{m=1}^{\infty}e^{-m\xi}\big[C_m\cos(m\eta)+D_m\sin(m\eta)\big],
\end{align}
which satisfies the condition $A_z^{\rm sc}\rightarrow0$ as $\xi\rightarrow\infty$. $C_m$ and $D_m$ are unknown expansion coefficients. To determine them, substitute \eqref{Azinc_elliptic_coords} and \eqref{seriesforaz} into \eqref{neumann-external-elliptic} and make use of the orthogonality of the trigonometric functions $\cos(m\eta)$ and $\sin(m\eta)$. We find
\begin{align}
\label{series-coefficients}
C_m=\frac{e^{m\xi_0}}{m\pi}\frac{b}{a}J\left(m,\frac{\rfil}{a},\frac{b}{a}\right),\quad D_m=0,\quad m=1,2,\dots.
\end{align}
where $J(m,x,y)$ is the convergent integral defined in (\ref{idefinition}). Using \eqref{series-coefficients} in conjunction with (\ref{ievaluated}) and some algebra, we can see that
\begin{align}
\label{Cm_coeff}
    C_m=\frac{e^{-m\xifil}-e^{m(2\xi_0-\xifil)}}{m}.
\end{align}
If we define
\begin{align}
w_1=e^{2\xi_0-\xifil-\xi},\quad w_2=e^{-\xifil-\xi},
\end{align}
then (\ref{seriesforaz}) becomes the convergent-series solution
\begin{align}
\label{seriesforaz1}
    A_z^{\rm sc}(\xi,\eta)=-\frac{\mu_0I}{2\pi}\bigg[-\sum_{m=1}^{\infty}\frac{w_1^m}{m}\cos(m\eta)+\sum_{m=1}^{\infty}\frac{w_2^m}{m}\cos(m\eta)\bigg], \quad \xi\geq\xi_0,
\end{align}
which somewhat resembles the series solution found in \cite{sten}. Finally, use of the following corollary of \eqref{identity3},
\begin{align}
\label{identity}
\sum_{m=1}^\infty\frac{w^m}{m}\cos(m\eta)=-\ln\sqrt{1-2w\cos{\eta}+w^2},\quad 0\le w<1, \quad -\pi<\eta\le\pi,
\end{align}
gives
\begin{align}\label{seriesforaz2}
    A_z^{\rm sc}(\xi,\eta)=-\frac{\mu_0I}{4\pi}\ln\frac{1-2w_1\cos\eta+w_1^2}{1-2w_2\cos\eta+w_2^2}.
\end{align}

\section*{Acknowledgments}
We thank Aggelos Dimopoulos for valuable help and useful discussions.

\bibliographystyle{siamplain}
\bibliography{mybib}
\end{document}